\documentclass[12pt,reqno]{amsart}

\usepackage[arrow,matrix,curve]{xy}

\usepackage[dvips]{graphicx} 

\usepackage{amssymb, latexsym, amsmath, amscd, array, hyperref
}

\theoremstyle{definition}

\numberwithin{equation}{section}

\newbox\GrossOneBox
\newbox\grossOneBox
\setbox\GrossOneBox\hbox{\raise-.4pt\hbox{${\rm O}\mskip-10.2mu\raise-2.7pt\hbox{$^1$}\mskip3mu$}}
\setbox\grossOneBox\hbox{\raise-4.5pt\hbox{\small$\mskip-1.3mu{\hbox{\rm\footnotesize O}}\mskip-10.0mu\raise-3.0pt\hbox{$^1$}\mskip1mu$}}

\author[P.B] {Piotr B\l{}aszczyk}\address{P. B\l{}aszczyk, Institute
of Mathematics, Pedagogical University of Cracow,
Poland}\email{pb@up.krakow.pl}

\author[A.B.]{Alexandre Borovik} \address{A. Borovik, School of
Mathematics, University of Manchester, Oxford Street, Manchester, M13
9PL, United Kingdom} \email{alexandre@borovik.net}

\author[V.K.]{Vladimir Kanovei} \address{V. Kanovei, IPPI, Moscow, and
MIIT, Moscow, Russia}\email{kanovei@googlemail.com}

\author[M.K.]{Mikhail G. Katz}\address{M. Katz, Department of
Mathematics, Bar Ilan University, Ramat Gan 52900 Israel}
\email{katzmik@macs.biu.ac.il}

\author[T.K.]{Taras Kudryk} \address{T. Kudryk, Department of
Mathematics, Lviv National University, Lviv, Ukraine}
\email{kudryk@mail.lviv.ua}

\author[S.K.]{Semen S. Kutateladze}\address{S. Kutateladze, Sobolev
Institute of Mathematics, Novosibirsk State University, Russia}
\email{sskut@math.nsc.ru}

\author [D.S.]  {David Sherry}\address{D. Sherry, Department of
Philosophy, Northern Arizona University, Flagstaff, AZ 86011, US}
\email{David.Sherry@nau.edu}

\begin{document}

\thispagestyle{empty}


\title [Analysis of a cultural icon: The case of Paul Halmos] {A
non-standard analysis of a cultural icon: The case of Paul Halmos}

\begin{abstract}
We examine Paul Halmos' comments on category theory, Dedekind cuts,
devil worship, logic, and Robinson's infinitesimals.  Halmos'
scepticism about category theory derives from his philosophical
position of naive set-theoretic realism.  In the words of an MAA
biography, Halmos thought that mathematics is ``certainty" and
``architecture" yet 20th century logic teaches us is that mathematics
is full of uncertainty or more precisely incompleteness.  If the term
\emph{architecture} meant to imply that mathematics is one great solid
castle, then modern logic tends to teach us the opposite lession,
namely that the castle is floating in midair.  Halmos' realism tends
to color his judgment of purely scientific aspects of logic and the
way it is practiced and applied.  He often expressed distaste for
nonstandard models, and made a sustained effort to eliminate
first-order logic, the logicians' concept of \emph{interpretation},
and the syntactic vs semantic distinction.  He felt that these were
\emph{vague}, and sought to replace them all by his \emph{polyadic
algebra}.  Halmos claimed that Robinson's framework is ``unnecessary''
but Henson and Keisler argue that Robinson's framework allows one to
dig deeper into set-theoretic resources than is common in Archimedean
mathematics.  This can potentially prove theorems not accessible by
standard methods, undermining Halmos' criticisms.


\end{abstract}

\maketitle

\tableofcontents

\section{Introduction}

Fifty years ago, the \emph{Pacific Journal of Mathematics} published a
pair of papers in the same issue, each containing a proof of a
conjecture in functional analysis known as the Smith--Halmos
conjecture.  The event had philosophical ramifications due to the fact
that one of these proofs involved methods that were not only unusual
for functional analysis but also challenged both historical thinking
about the evolution of analysis and foundational thinking in
mathematics.  The present article explores these and related issues.

Paul Halmos was a 20th century expert in functional analysis.  His
textbooks on measure theory, Hilbert spaces, and finite dimensional
vector spaces are well written, still relevant, and highly praised.%
\footnote{See
\url{http://www.maa.org/news/paul-halmos-a-life-in-mathematics}}

Following the Aronszajn--Smith proof of the existence of invariant
subspaces for compact operators \cite{AS}, Smith and Halmos
conjectured that the same should be true for more general classes of
operators, such as operators with a compact square.  A proof in the
more general case of polynomially compact operators in \cite{BR}
(exploiting Robinson's infinitesimals) was a notable event in
functional analysis.  Simultaneously the same journal published an
infinitesimal-free proof \cite{Ha66}.

In 1991, J. Dauben interviewed the distinguished model theorist
C.~C.~Chang about the Bernstein--Robinson paper.  Even a quarter
century later (and after \cite{Lo73} superseded the 1966 results)
Chang still seemed a bit sore about Bernstein and Robinson not getting
enough credit, for he insisted that
\begin{quote}
once you know something is true, it is easier to find other proofs.
Major credit must go to Robinson.%
\footnote{Chang's reference to Robinson is certainly shorthand for
Bern\-stein--Robinson.}
(Chang quoted in \cite[p.~327]{Da95})
\end{quote}
Robinson himself supports Chang's reading:
\begin{quote}
As for the Halmos standard `translation', it was all very nice, but
the NSA (i.e., \emph{nonstandard analysis}) proof was quite natural,
while the standard proof required an argument that would not have been
so easy to spot without first seeing the NSA version.  (Reported by
Moshe Machover, private communication)
\end{quote}
In a course at Hebrew University in the late 1960s, Robinson said:
\begin{quote}
Halmos was proud of his proof but in the end all he did was rewrite
our proof in a language he was educated in.  (Reported by Shmuel
Dahari, private communication).
\end{quote}
Halmos himself essentially agreed with this sentiment when he wrote
that the purpose of his paper was 
\begin{quote}
to show that by appropriate \emph{small modifications}[,] the
Bernstein--Robinson proof can be converted \ldots{} into one that is
expressible in the standard framework of classical analysis.
\cite[p.\;433]{Ha66} (emphasis added) 
\end{quote}
Further details can be found in Section~\ref{s33}.

Subsequently Halmos expressed reservations about Robinson's framework,
and described researchers working in the framework as \emph{converts}
(see Section~\ref{s61}).

What philosophical outlook shaped Halmos' attitude toward Robinson's
framework, and prompted his critical remarks concerning fellow
experts?  Following \cite{Je04}, we provide an analysis that can
hardly be described as \emph{standard} of a little-known aspect of a
mathematical \emph{cultural icon}.

\section{Paraphrase}
\label{s33}

The invariant subspace conjecture of Smith and Halmos was first proved
by Bernstein and Robinson, and published in the \emph{Pacific Journal
of Mathematics} (\emph{PJM}).  A number of scholars would have been
more comfortable had Halmos' infinitesimal-free paraphrase of the
proof in \cite{BR} (for which Halmos was apparently the referee)
appeared in the \emph{next} issue of the \emph{PJM} rather than being
published simulataneously in the same issue as \cite{Ha66}.

The Bernstein--Robinson proof is presented in detail in \cite{Da77}.

Halmos claimed two decades later that he received the manuscript by
Bernstein and Robinson ``early in 1966'' in \cite[p.~320]{Ha85}, but
that date is certainly incorrect.  Dauben documents a letter from
Halmos to Robinson acknowledging receipt of the manuscript, and dated
19~june 1964 (see \cite[p.~328, note 66]{Da95}).  Thus Halmos was in
possession of the Bernstein--Robinson manuscript even \emph{prior to}
its submission for publication on 5~july 1964.

Several specialists have privately testified that Halmos was most
likely the referee for the Bernstein--Robinson paper.%
\footnote{See also a related discussion at
\url{http://mathoverflow.net/questions/225455}.}
A slightly delayed publication of Halmos' paraphrase (say, in the
following issue of the \emph{PJM}) would have avoided the effect of
weakening the Bernstein--Robinson priority claim on the result, and
may have constituted a more appropriate use of publication timetables.
There have been several cases of scholars affected by the
marginalisation campaign against Robinson's framework who ended up
suffering in terms of employment as a result, indicating that such
issues are not purely \emph{academic}.

Here by ``Robinson's framework'' we mean Robinson's rigorous
justification of Leibnizian infinitesimal procedures in the framework
of modern mathematics (viz., the Zermelo--Fraenkel set theory with the
axiom of choice), as developed in \cite{Ro61} and \cite{Ro66}.
Robinson exploited the theory of types in presenting his framework.
Alternative presentations involve ultraproduct constructions; see
e.g., \cite{Lu62}.

\section{Indispensability argument of Henson and Keisler}
\label{s333}

Halmos explicitly referred to his own paper as a ``translation'' (of
the Bernstein--Robinson proof).  However he did not think of it as an
\emph{awkward} translation, and on the contrary used it to justify his
claim in \cite{Ha85} that NSA is unnecessary because it can always be
translated.  The following year, Henson and Keisler published a paper
\cite{HK} that was a reaction to a widespread belief at the time that
Robinson's framework is unnecessary, and in particular provided a
rebuttal of Halmos' claims.

\subsection{Second-order Arithmetic}

Henson and Keisler point out that a nonstandard extension of second
order Arithmetic is not a conservative extension of second-order
Arithmetic, but is rather closely related to third-order theory.  This
is because, roughly, nonstandard arguments often rely on
\emph{saturation} techniques that typically involve third-order
theory.  They go on to argue against the type of fallacy contained in
Halmos' position that Robinson's framework is unnecessary.  The gist
of their argument is that since most mathematics takes place at
second-order level, there may well be nonstandard proofs whose
standard translations, while theoretically possible, may well be
humanly incomprehensible.  They conclude as follows:
\begin{quote}
This shows that in principle there are theorems which can be proved
with nonstandard analysis but cannot be proved by the usual standard
methods. The problem of finding a specific and mathematically natural
example of such a theorem remains open.  \cite[p.~377]{HK}
\end{quote}
In this spirit, \cite{TV} use the language of Robinson's framework in
order to avoid a large number of iterative arguments to manage a large
hierarchy of parameters.  Ultraproducts form a bridge between discrete
and continuous analysis \cite{Go97}.

\subsection{Rebuttal of Halmos' claims}
\label{s42}

Halmos formulated a pair of claims concerning Robinson's framework,
which are closely related but perhaps not identical:
\begin{enumerate}
\item
Robinson had a language and not an idea.
\item
Robinson introduced a special tool, too special, and other tools can
do everything it can, so it's all a matter of taste.
\end{enumerate}
In Halmos' own words:
\begin{quote}
If they had done it in Telegu [sic]%
\footnote{The correct spelling is \emph{Telugu}.}
instead, I would have found their paper even more difficult to decode,
but the extra difficulty would have been one of degree, not of kind.
\cite[p.~204]{Ha85}
\end{quote}

Even though Halmos calls it a ``language" in (1) and a ``tool" in (2),
the underlying claim is essentially the same: just as you can express
your mathematics in English, French, or Telugu and it does not make
any difference, so similarly you can do your mathematics in
traditional set-ups or in a Robinsonian logical contraption.

The rebuttal is the same in both cases, and was already provided by
the Henson--Keisler argument and the example of Tao's work, as
discussed in Section~\ref{s333}.

Today, Robinson's framework is neither a language, idea, or tool, but
rather is a branch of modern mathematics with its own domain, set of
tools, collection of key results, and numerous applications.

\section{Dedekind cuts and category theory}
\label{s3}

The following comment by Halmos needs to be addressed:
\begin{quote}
Here is a somewhat unfair analogy: Dedekind cuts. It's unfair because
it's even more narrowly focused, but perhaps it will suggest what I
mean. No, we don't have to learn it (Dedekind cuts or non-standard
analysis): it's a special tool, too special, and other tools can do
everything it does. It's all a matter of taste.  \cite[p.~204]{Ha85}
\end{quote}
Halmos seems to view both Dedekind cuts and category theory with
disfavor.  On the other hand, one who doesn't favor cuts should
apparently favor category theory, since excising cuts would make the
real line a \emph{category}, i.e., something without a strict
set-theoretic definition.

\subsection{Category theory viewed by some}

Halmos' attitude to Robinson's framework is somewhat comparable to
Halmos' attitude to category theory, at the expense of which he also
made disparaging remarks:
\begin{quote}
A microscopic examination of such similarities might lead to category
theory, a subject that is viewed by some with the same kind of
suspicion as logic, but not to the same extent.  \cite[p.~205]{Ha85}
\end{quote}
In his essay ``Applied mathematics is bad mathematics,'' Halmos
claimed that when applied mathematicians describe category theory as
``abstract nonsense,'' they mean it \cite[p.~15]{Ha81}, but provided
no evidence to substantiate his claim that applied mathematicians feel
this way, or that such sentiments are due to anyone but himself.

Halmos sought to identify categories with universal algebras, thus
reducing category theory to set theory in \cite{Ha81b}.

Category theory is today one of the fastest growing industries, with
avid advocates like David Kazhdan.  Halmos might have pigeon-holed
Kazhdan a ``convert'' as well (see Section~\ref{s61}), but it wouldn't
have helped Halmos' reputation.

\subsection{Bridge between discrete and continuous}

Robinson's framework is a fruitful modern research area that has
attracted many researchers, as noted in Section~\ref{s42}.  Halmos
predicted that
\begin{quote}
in the foreseeable future \ldots{} \emph{discrete} mathematics will be
an increasingly useful tool in the attempt to understand the world,
and \ldots{} analysis will therefore play a proportionally smaller
role.  \cite[p.~19]{Ha81} (emphasis added)
\end{quote}
What Halmos may not have anticipated is that, in fact, the
ultraproducts form a bridge between discrete and continuous analysis
as mentioned above.

\section{Halmos and logic}

The algebraic approach to logic has a long history starting with
Boole, continuing with Peirce and Schr\"oder, and reaching a high
point with the L\"owenheim--Skolem theorem.  Subsequently it went out
of fashion to a certain extent, but the work of Tarski on Boolean
algebras with operators eventually led to his cylindric algebras,
i.e., Bolean algebras with \emph{quantifiers} as the added operators.
The Tarski school has proved a number of difficult, and perhaps even
\emph{deep}, results about this class of algebras.

\subsection{From cylindric to polyadic algebras}
\label{s61b}

Halmos became interested in this topic, as he discusses in his book
\cite{Ha85}, where one finds some remarks on polyadic vs cylindric
algebras; see also \cite{Ha00}.  Whether or not there are any
contributions of substance by Halmos to logic proper is a delicate
question.  His polyadic formalism differs from the cylindric
counterpart, but the theory in his book is a straightforward
translation of first order logic, thus not \emph{deep} by any means.
Neither polyadic nor cylindric algebras made a major contribution to
logic and its applications, and are of marginal interest today.

In later work on probability, the algebraic formalism was dropped in
favor of working within first order logic.  Halmos' translation of the
completeness theorem, i.e., his representation theorem, is rather
complicated.  Thus, Fenstad gave a simplified presentation and used
this work to give a rather general representation theorem for logical
probabilities in \cite{Fe67}.

Halmos' feelings about logic in general and Robinson's framework in
particular are neatly summarized in a limerick dating from 1957, and
republished on page~216 in his book:

\begin{verse}
If you think that your paper is vacuous,\\
Use the first-order functional calculus.\\
It then becomes logic,\\
And, as if by magic,\\
The obvious is hailed as miraculous.
\end{verse}

It has to be admitted that Halmos and Bishop had something in common,
namely literary talent (see Section~\ref{s72}).  The limerick aptly
summarizes the import of Halmos' own contribution to logic.

\subsection{Quixotic battle against formal logic}

A passage in Halmos' book reproduced in his article ``An Autobiography
of Polyadic Algebras" is part of his attack on \emph{formal logic} (as
opposed to \emph{symbolic logic} favored by Halmos), and runs as
follows:
\begin{quote}
When I asked a logician what a variable was, I was told that it was
just a `letter' or a `symbol'. Those words do not belong to the
vocabulary of mathematics; I found the explanation that used them
unhelpful--\emph{vague}.  When I asked what `interpretation' meant, I
was answered in bewildering detail (set, correspondence, substitution,
satisfied formulas).  In comparison with the truth that I learned
later (homomorphism), the answer seemed to me unhelpful--forced, ad
hoc. It was a thrill to learn the truth--to begin to see that formal
logic might be just a flat photograph of some \emph{solid}
mathematics--it was a thrill and a challenge.  \cite[p.~208]{Ha85},
\cite[p.~385-386]{Ha00} (emphasis added)
\end{quote}
Some issues need to be clarified in connection with this passage:
\begin{enumerate}
\item
What is Halmos' problem with formal logic exactly?
\item
What is wrong with the term \emph{interpretation}?
\item
In what way does replacing the term \emph{interpretation} by the term
\emph{homomorphism} help?
\item
What is \emph{unsolid} about formal logic?
\end{enumerate}
Exploring these questions may help understand Halmos' 36 year battle
(1964--2000) against anything nonstandard.%
\footnote{There is yet another dig against non-standard models in his
2000 article cited above, one of the last ones he wrote.}
What Halmos seems to be reacting against is a distinction taken for
granted in modern logic, namely that between syntax and semantics.
Roughly, this means that one can have a \emph{theory} at the syntactic
level which does not \emph{mean} anything until one \emph{interprets}
it in a specific model to get meaning (\emph{semantics}).  This view
presupposes a possibility of having \emph{distinct} models for the
same \emph{theory}.

\subsection{Mathematics as one great thing}

Halmos' position against such dualities appears to stem from a naive
set-theoretic realism (already on display in his opposition to
category theory; see Section~\ref{s3}).  Halmos seems opposed to the
idea that there are distinct levels of things in mathematics: you can
have a \emph{theory} of a distinct level of mathematical \emph{Sein}
than an \emph{interpretation} thereof.  Halmos apparently prefers to
see all mathematics as made of the same cloth:
\begin{quote}
I see mathematics, the part of human knowledge that I call
mathematics, as one thing--one great, glorious thing.  (Halmos quoted
in \cite[p.~234]{Al82})
\end{quote}
Now the `one great thing' comment suggests that all mathematical
objects are sets, and sets differ in degree of complexity but they do
not differ in kind.

In this sense, \emph{homomorphisms} are more \emph{solid} than
\emph{interpretations}, in that talking about homomorphisms implies
that the domain and the range are of the same kind, thereby escaping
the duality of theory/interpretation that seems to threaten the
solidity of naive set-theoretic realism.  Perhaps Halmos' polyadic
algebras could be understood as an attempt to undo \emph{formal logic}
with its threatening dualities and inherent possibilities of unsolid
(a.k.a., nonstandard) models.  A related point was made by G. Lolli,
in the context of an analysis of Halmos' views, in the following
terms:
\begin{quote}
\ldots{} the deep reason for the opposition, depreciation and
misunderstandings concerning logic among mathematicians lies in their
inability or unwillingness to accept the binomium
language-metalanguage as a mathematical tool; they don't even seem
capable of understanding its sense. This could be due to their habit
of talking in an informal quasi-natural language, where metalanguage
is flattened on the language itself, or the languages are absorbed in
the metalanguage, a habit legitimated and reinforced by the
set-theoretical framework.  \cite{Lo08}
\end{quote}
Having identified the set-theoretic source of the problem, Lolli
concludes: 
\begin{quote}
They should know however, as everybody is now aware, that this very
identification is the source of dangerous circularities. Only the
conceptual distinction, at least in principle, of language and
metalanguage avoids the paradoxes. (ibid.)
\end{quote}

\section{A rhetorical analysis}

In addition to scientific arguments, Halmos resorted on several
occasions to excesses of language aimed at marginalizing Robinson's
framework, as we document in this section.

\subsection{Halmos on types of worship}
\label{s61}

Halmos may have been a leading expert in his field, but so was Edward
Nelson (see e.g., \cite{Ne67}), and so is Peter Loeb (see e.g.,
\cite{Lo75}).  Halmos had the following to say about their relation to
Abraham Robinson's framework:
\begin{quote}
\ldots for some \emph{converts} (such as Pete Loeb and Ed Nelson),
it's a religion, \ldots{} For some others, who are against it (for
instance Errett Bishop), it's an equally emotional issue--they regard
it as \emph{devil worship}.  \cite[p.~204]{Ha85} (emphasis added)
\end{quote}

Halmos' description of both Nelson and Loeb as ``converts" in the
comment quoted above raises questions of motivation behind applying
this kind of epithet to fellow leading mathematicians, or for that
matter of invoking Errett Bishop on ``devil worship,'' remarks that
are dangerously close to the category of expletives.  In point of fact
Bishop never used such a term in reference to Robinson's
infinitesimals (see more on \emph{devil worship} in
Section~\ref{s72}).  Halmos sought to create the impression of a
balanced presentation of \emph{both sides} of the controversy by
mentioning \emph{both} Nelson \emph{and} Bishop, but in fact both of
his \emph{sides} serve only as a vehicle for an attempt to demonize
Robinson's framework.

\subsection{Errett Bishop}
\label{s72}

Halmos' remarks concerning \emph{devil worship} in Section~\ref{s61}
deserve closer scrutiny.  Bishop's verse on the \emph{neat devil} that
is classical mathematics, from his essay ``Schizophrenia in
contemporary mathematics,'' run as follows:
\begin{verse}
The devil is very neat. It is his pride \\ To keep his house in
order. Every bit \\ Of trivia has its place. He takes great pains \\
To see that nothing ever does not fit.  \\ And yet his guests are
queasy.  All their food, \\ Served with a flair and pleasant to the
eye, \\ Goes through like sawdust.  Pity the perfect host!  \\ The
devil thinks and thinks and he cannot cry.  
\end{verse}
(See \cite[p.\;14]{Bi85}.)  For additional details on Bishop's antics
see \cite{KK11d}, \cite{KK12a}, \cite{Ka15}.  The ``Schizophrenia''
essay says not a word about Robinson's framework, and all the
\emph{devil} material (verse or prose) targets classical mathematics
\emph{as a whole}, including Halmos' favorate subjects such as
invariant subspaces.  Bishop's poem was published earlier but composed
later than his teacher%
\footnote{Apparently in more than one area}
Halmos' limerick; see Section~\ref{s61b}.  Halmos' claim that Bishop
regarded Robinson's framework as devil worship appears to be merely a
smear-by-proxy attack on Robinson.  It is certainly possible that
Bishop may have made private remarks along these lines to Halmos, who
was after all his advisor.  Still, Halmos' purported quote of Bishop
cited in Section~\ref{s61} is taken out of context.

We are not sure whether there is an official philosophical term for
such a rhetorical technique, but at any rate it is not the unique
occurrence of such a technique in Halmos.  He did something similar
with regard to category theory, while positioning himself safely
behind the broad backs of unnamed applied mathematicians; see
Section~\ref{s3}.

\subsection{Underworld}

What would be the point of using mocking epithets like ``dredged up
from the underworld,'' as Halmos did in his 1990 article, in
describing Robinson's accomplishment with regard to infinitesimals:
\begin{quote}
The modern theory of nonstandard analysis \emph{dredged} the forbidden
concepts \emph{up from the underworld} and is trying to reinstate them
at the right side of Cauchy's throne.  \cite[p.~569]{Ha90} (emphasis
added)
\end{quote}
Halmos may have been more moderate in his language than Connes who
used some objectionable vitriol in referring to Robinson's framework
(see \cite{KKM}, \cite{KL}), but in the end Halmos' attitude is
comparable to Connes', that other leading expert.  In fact, in his
book Halmos broadened his criticism of Robinson to a broader criticism
of logic:
\begin{quote}
The logician's attention to the nuts and bolts of mathematics, to the
symbols and words ($0$ and $+$ and ``or'' and ``and''), to their order
($\forall\exists$ or $\exists\forall$), and to their grouping
(parentheses) can strike \emph{the} mathematician as pettifogging
\ldots{} \cite[p.~205]{Ha85} [emphasis added]
\end{quote}
The definite article attached to ``mathematician'' is the issue here,
for it presupposes that there is just one thing that counts as being a
mathematician.  `Some' would make it more accurate, but significantly
blunt the force of the remark.

Here Halmos is apparently alluding to Robinson's approach to
infinitesimals via the theory of types, with its reliance on the
``nuts and bolts'' of logic.  If Halmos wished to publish an
evaluation of Robinson's framework, he could have been expected to
have done enough research to discover a more elementary analytical
approach.  This is the ultrapower approach, already exploited in
\cite{He48} and popularized by Luxemburg in the CalTech Lecture Notes
and e.g., in \cite{Lu62}, namely over two decades prior to the
publication of Halmos' book.

The sweeping and sarcastic critique Halmos presents fails to inform
the reader that there does exist an accessible analytical approach to
infinitesimals \cite{Li88}.  The existence of such an approach makes
much of Halmos' vitriol rather misplaced.  There might exist more
abstract approaches that he does not appreciate, but the same can be
said about many fields in mathematics.  There are certainly textbooks
in, for example, differential geometry that are more accessible than
other textbooks in differential geometry.  The existence of the more
abstract textbooks generally does not lead sceptical scholars to speak
of differential geometry as being ``dredged up from the underworld.''

\section{Conclusion}

Robinson's characterisation of Bishop's ``attempt to describe the
philosophical and historical background of [the] remarkable endeavor''
of the constructive approach to mathematics, as ``more vigorous than
accurate'' \cite[p.~921]{Ro68} applies equally well to Halmos' take on
logical issues, conditioned by his naive set-theoretic realism.  Such
a philosophical \emph{parti pris} led Halmos to reject not merely
Robinson's infinitesimals but also broad swaths of standard techniques
and applications, ranging from a modern logical toolkit like
first-order logic to applied mathematics.  Halmos' attempted reform of
logic is a radical project that bears similarity to his student Errett
Bishop's even more radical opposition to classical mathematics as a
whole, as analyzed elsewhere.

\end{document}